\documentclass[11pt]{amsart}
\usepackage{graphics}
\usepackage{amsmath}
\usepackage{amsfonts}
\usepackage{amsthm}
\usepackage{latexsym}
\usepackage{amssymb}
\newcommand{\g}{{\mathfrak g}}
\newcommand{\h}{{\mathfrak h}}
\newcommand{\n}{{\mathfrak n}}

\newcommand{\C}{{\bf C}}
\newcommand{\R}{{\bf R}}

\newcommand{\Z}{{\bf Z}}

\newtheorem{theorem}{\bf Theorem}
\newtheorem{proposition}{\bf Proposition}
\newtheorem{lemma}{\bf Lemma}

\newtheorem{remark}{\bf Remark}
\newtheorem{example}{\bf Example}
\newtheorem*{conjecture}{\bf Conjecture}
\begin{document}
%

\title[Non-left-invariant complex structures]
{Small deformations and non-left invariant complex structures
on a compact solvmanifold
}
\author{Keizo Hasegawa
}


\date{}

\maketitle

\begin{abstract}
We observed in our previous paper that all the complex structures on
four-dimensional compact solvmanifolds, including tori, are left-invariant. 
In this paper we will give an example of a six-dimensional compact
solvmanifold which admits a continuous family of non-left-invariant complex
structures. Furthermore, we will make a complete classification of three-dimensional compact
homogeneous complex solvmanifolds; and determine which of them
admit pseudo-K\"ahler structures.
\end{abstract}



\section{Introduction}

A homogeneous manifold $M$ is a differentiable manifold on which a real
Lie group $G$ acts transitively. In the case where $M$ is a complex manifold,
we usually assume that the group action is holomorphic, and $M$ is called a
{\em homogeneous complex manifold}. It should be noted \cite{W2}
that any compact homogeneous complex manifold can be written as
$H \backslash G$, where $G$ is a complex Lie group and $H$ is a closed
complex subgroup of $G$.
A homogeneous complex structure on a Lie group $G$ (considered as
a homogeneous manifold) is nothing but a left-invariant complex structure on $G$;
and it is a complex Lie group if and only if it is both left and right-invariant.
There are already extensive studies on left-invariant complex structures on
Lie groups (\cite{OR}, \cite{S}, \cite{Sn}): for instance, the classification of all
homogeneous complex surfaces is known, which includes all two-dimensional compact
homogeneous complex manifolds and all left-invariant
complex structures on four-dimensional simply connected Lie groups \cite{OR}.
Any left-invariant complex structure on a Lie group $G$ defines a canonical complex
structure on its quotient $M = \Gamma \backslash G$, where $\Gamma$ is
a discrete subgroup of $G$.
We call such a complex structure a {\em left-invariant complex structure} on $M$.
Remark that unless  the canonical right action of $G$ on $M$ is holomorphic,
$M$ may not be a homogeneous complex manifold. 
\smallskip

In this paper we call a  compact homogeneous manifold of
solvable (nilpotent) Lie group a {\em solvmanifold} ({\em nilmanifold});
and a homogeneous complex solvmanifold (nilmanifold)
a {\em complex solvmanifold} ({\em complex nilmanifold}).
In our previous paper \cite{H1}, we showed that all the complex structures on
a four-dimensional solvmanifold $M$ are left-invariant: that is, expressing 
$M$ as $\Gamma \backslash G$ (up to finite covering), where $G$
is a four-dimensional simply connected solvable Lie group and $\Gamma$ is a
lattice of $G$, they are all induced from some left-invariant complex structures on $G$.
To be more precise, we showed:

\begin{theorem}[\cite{H1}]
A complex surface is diffeomorphic to a four-dimensional solvmanifold
if and only if it is one of the following surfaces:
Complex torus, Hyperelliptic surface, Inoue Surface of type $S^0$,
Primary Kodaira surface, Secondary Kodaira surface, Inoue Surface of type
$S^{\pm}$. Furthermore, every complex structure on each of these complex surfaces
(considered as solvmanifolds) is left-invariant.
\end{theorem}

A natural question then arises whether the last assertion in the theorem also holds for
higher dimension. We will show in Section 5 that there
exists an example of a
six-dimensional solvmanifold which admits a continuous family of non-left-invariant
complex structures (see Theorem 4). In fact, in the paper \cite{N} Nakamura constructed
small deformations of a three-dimensional complex solvmanifold
$\Gamma \backslash G$,  where $G$ is a complex solvable Lie group of dimension 3 and
$\Gamma$ is a lattice (uniform discrete subgroup) of $G$; and showed in particular that
there exists a continuous family of complex structures whose universal coverings are
not Stein (as noted in the paper, this construction is actually due to Kodaira). 
Therefore, in order to show that there exist non-left-invariant complex structures on a
six-dimensional  solvmanifold, it is sufficient to show that all the left-invariant
complex structures on $G$ are biholomorphic to ${\bf C}^3$ (see Theorem 3).
Note that this result implies that neither complex-homogeneity nor left-invariance of
complex structure is preserved under small deformations. 

We have some conjectures relating to small deformations and left-invariant complex
structures on  solvmanifolds.

\begin{conjecture}
(i) All the left-invariant complex structures on even-dimensional simply connected
unimodular solvable Lie groups (nilpotent Lie groups)
are Stein (biholomorphic to $\C^n$ respectively);
(ii) Small deformations of left-invariant complex structures on
even-dimensional nilmanifolds are all left-invariant.
\end{conjecture}

It should be noted that a simply connected complex solvable Lie group of
dimension $n$ is biholomorphic to $\C^n$;
and that small deformations of a complex torus are all left-invariant. 
Remark also that the conjectures (i) and (ii) hold for dimension 4 (see \cite{H1}, \cite{OR}). 
\smallskip

In the paper \cite{N} Nakamura has classified three-dimensional  complex
solvmanifolds $M$ into four classes: (1) abelian type with $h^1 = 3$,
(2) nilpotent type with $h^1 = 2$, (3a) non-nilpotent type with $h^1 = 1$,
and (3b) non-nilpotent type with $h^1 = 3$,
where $h^1 = {\rm dim}\, H^1(M, {\mathcal O})$. In Section 3,
determining all lattices of three-dimensional unimodular complex solvable Lie groups
we will complete the classification of three-dimensional complex solvmanifolds.
\smallskip

In the paper \cite{Y1} Yamada gave the first example of a complex solvmanifold
which admits a pseudo-K\"ahler structure. In Section 4 we will show that a three-dimensional
complex solvmanifold admits a pseudo-K\"ahler structure if and
only if it is of type (1) or (3b) (see Theorem 2).


\section{Preliminaries}

Let $G$ be a connected, simply connected Lie group of dimension $2m$, and
$\g$ the Lie algebra of $G$. We consider a left-invariant almost complex structure
$J$ on $G$ as a linear automorphism of $\g$, that is, $J \in {\rm GL} (\g, \R)$
such that $J^2 = - I$. As is well known $J$ is integrable (that is, it defines
a left-invariant complex structure on $G$) if and only if the Nijenhuis tensor $N_J$
on $\mathfrak g$ vanishes identically, where $N_J$ is defined by
$$N_J(X,Y)=[JX,JY]-J[JX,Y]-J[X,JY]-[X,Y]
$$
for $X, Y \in \g$.

Let $\g_\C =\g \otimes \C$ be the complexification of $\g$.
We will reformulate the integrability condition of $J$ in terms of complex
subalgebras of $\g_\C$. For an almost complex structure $J$ on $\g$,
let $\h_J$ be the complex subspace of $\g_\C = \g \oplus \sqrt{-1}\g$ generated by
$X + \sqrt{-1} JX, X \in \g$, that is,
$$\h_J = \{X + \sqrt{-1} JX| X \in \g \}_\C\;.
$$
Then, we see that $J$ is integrable if and only if $\h_J$ is a complex Lie subalgebra
of $\g_\C$ such that $\g_\C = \h_J \oplus \overline{\h_J}$. On the other hand, we have
the following lemma.

\begin{lemma}
Let $V$ be a real vector space of dimension $2m$.
Then, for a complex subspace $W$ of $V \otimes \C$ such that
$V \otimes \C = W \oplus \overline{W}$, there is a unique $J_W \in {\rm GL}(V, \R),
{J_W}^2 = -I$ such that
$$W = \{X + \sqrt{-1} J_W X| X \in V \}_\C.$$
\end{lemma}

It follows from this Lemma and the above argument that there exists
one to one correspondence between complex
(integrable almost complex) structures on $\g$ and complex Lie subalgebras
$\h$ such that $\g_\C = \h \oplus \overline{\h}$. The correspondence is given by
$J \rightarrow \h_J$ and $\h \rightarrow J_\h$. 

We now suppose that $J$ is a complex structure
on $\g$ with its associated complex Lie algebra $\h_J$. Then the complex Lie subgroup
$H_J$ of $G_\C$ corresponding to $\h_J$ is closed, simply connected,
and $H_J \backslash G_\C$ is biholomorphic to $\C^m$. The canonical inclusion
$\g \hookrightarrow \g_\C$ induces an inclusion $G \hookrightarrow G_\C$,
and $\Gamma = G \cap H_J$ is a discrete subgroup of $G$.
We have the following canonical map $g = i \circ \pi$:
$$G \stackrel{\pi}{\rightarrow} \Gamma \backslash G \stackrel{i}{\hookrightarrow}
H_J \backslash G_\C,
$$
where $\pi$ is a covering map, and $i$ is an inclusion. We can see that
${\rm Im}\; g$ is an open subset $U$ of $\C^m$, and the complex structure
$J$ on $G$ is the one induced from $U \subset \C^m$ by $g$.
It should be noted that if $G$ is a complex Lie group, we have
$\Gamma = G \cap H_J = \{1\}$, and $g$ is a biholomorphic map onto $\C^m$.
For the details of the above argument we refer to the paper \cite{Sn}.


\section{Three-dimensional unimodular complex solvable Lie groups}

A complex solvmanifold can be written as
$\Gamma \backslash G$, where $G$ is a simply connected, unimodular complex
solvable Lie group and $\Gamma$ is a lattice (uniform discrete subgroup) of $G$ \cite{BO}.
In particular, the Lie algebra $\g$ of $G$ must be unimodular; that is, the trace of ad (X) is
$0$ for every $X$ of $\g$. It is easy to classify
all unimodular complex solvable Lie algebras of dimension 3.
They are divided into three classes: (1) abelian type, (2) nilpotent type,
(3) non-nilpotent type.

In the following list, we express the solvable Lie algebra $\g$ as having
a basis $\{X, Y, Z\}$  with the bracket multiplication specified for each type:

\begin{list}{}{\topsep=5pt \leftmargin=15pt \itemindent=5pt \parsep=0pt \itemsep=3pt}
\item[ (1)] Abelian Type: \par
$[X, Y] = [Y, Z] = [X, Z] = 0$.
\item[ (2)] Nilpotent Type: \par
$[X, Y] = Z,\; [X, Z] = [Y, Z] = 0$.
\item[ (3)] Non-Nilpotent Type: \par
$[X, Y] = -Y,\; [X, Z] = Z,\; [Y, Z] = 0$.
\end{list}

For each of their corresponding simply connected solvable Lie groups $G$,
we will determine all lattices $\Gamma$:

\begin{list}{}{\topsep=5pt \leftmargin=5pt \itemindent=5pt \parsep=0pt \itemsep=5pt}

\item[ (1)] Abelian Type:\; $G = \C^3$ \par
A lattice $\Gamma$ of $G$ is generated by a basis of $\C^3$ as a vector space
over $\R$, and $\Gamma \backslash G$ is a complex torus.

\item[ (2)] Nilpotent Type:\; $G = \C^2 \rtimes \C$ with the action $\phi$ defined by
$$\phi(x)(y, z) = (y, z + xy),
$$
or in the matrix form,
$$G = \Bigg\{ \left(
\begin{array}[c]{ccc}
1 & x & z\\
0 & 1 & y\\
0 & 0 & 1
\end{array}
\right) \rule[-7mm]{0.25mm}{16mm}\;  x, y, z \in \C \Bigg\}.
$$
A lattice $\Gamma$ of $G$ can be written as 
$$\Gamma = \Delta \rtimes \Lambda,
$$
where $\Delta$ is a lattice of $\C^2$ and $\Lambda$ is a lattice of $\C$.
Since an automorphism $f \in {\rm Aut}(\C)$ defined by
$f(x) = \alpha x, \alpha \not= 0$ can be extended to
an automorphism $F \in {\rm Aut}(G)$ defined by $F(x, y, z) = (\alpha x, \alpha^{-1} y, z)$,
we can assume that $\Lambda$ is generated by $1$ and
$\lambda\; (\lambda \notin \R)$ over
$\Z$. Since $\Delta$ is preserved by $\phi(1)$ and $\phi(\lambda)$, we see that
$\Delta$ is generated by
$(\alpha_1, \beta_1), (\alpha_2, \beta_2), (0, \alpha_1),$ $(0, \alpha_2)$ over $\Z$,
where $\beta_1$ and $\beta_2$ are arbitrary complex numbers, and $\alpha_1$ and
$\alpha_2$ are linearly independent over $\R$ such
that $(\alpha_1, \alpha_2)$ is an eigenvector of some 
$A \in {\rm GL}(2, \Z)$ with the eigenvalue $\lambda$. Conversely, for any
$A \in {\rm GL}(2, \Z)$ with non-real eigenvalue $\lambda$, we can define a lattice
$\Gamma$ of $G$.

\begin{example}
A standard lattice $\Gamma$ of $G$ with $x, y, z \in \Z[\sqrt{-1}]$ is obtained by
putting $\lambda = \sqrt{-1}, \alpha_1 = \alpha_2 = 0, \beta_1 = 1, \beta_2 = \sqrt{-1}$,
and $\Gamma \backslash G$ is an {\em Iwasawa manifold}.
\end{example}

\item[ (3)] Non-Nilpotent Type:\; $G = \C^2 \rtimes \C$ with the action $\phi$
defined by
$$\phi(x)(y, z) = (e^{x} y, e^{-x} z),
$$
or in the matrix form,
$$G = \Bigg\{ \left(
\begin{array}[c]{cccc}
e^{x} & 0 & 0 & y\\
0 & e^{-x} & 0 & z\\
0 & 0 & 1 & x\\
0 & 0 & 0 & 1
\end{array}
\right) \rule[-7mm]{0.25mm}{16mm}\;  x, y, z \in \C \Bigg\}.
$$

A lattice $\Gamma$ of $G$ can be written as $\Gamma = \Delta \rtimes \Lambda$,
where $\Delta$ is a lattice of $\C^2$, and $\Lambda$ is a lattice of $\C$
which is generated by $\lambda$ and $\mu$ over $\Z$. 
Since $\Delta$ is preserved by $\phi(\lambda)$ and $\phi(\mu)$,
we see that $\Delta$ is generated by $(\alpha_i, \beta_i), i = 1, 2, 3, 4$ over $\Z$
such that 
$$\textstyle \gamma^{-1} \alpha_i = \sum_{j=1}^{4} a_{i j} \alpha_j,\; 
\gamma \beta_i = \sum_{j=1}^{4} a_{i j} \beta_j,
$$
$$\textstyle \delta^{-1} \alpha_i = \sum_{j=1}^{4} b_{i j} \alpha_j,\; 
\delta \beta_i = \sum_{j=1}^{4} b_{i j} \beta_j,
$$
where $\gamma = e^{\lambda}, \delta = e^{\mu}$, and
$A = (a_{i j}), B = (b_{i j}) \in {\rm SL}(4, \Z)$ are semi-simple and
mutually commutative.
In other word, we have simultaneous eigenvectors 
$\alpha = (\alpha_1, \alpha_2, \alpha_3, \alpha_4),\\
\beta = (\beta_1, \beta_2, \beta_3, \beta_4) \in \C^4$ 
of $A$ and $B$ with eigenvalues $\gamma^{-1}, \gamma$ and $\delta^{-1}, \delta$
respectively.
Conversely, for any mutually commutative, semi-simple matrices
$A, B \in {\rm SL}(4, \Z)$ with eigenvalues $\gamma^{-1}, \gamma$ and
$\delta^{-1}, \delta$ respectively, take simultaneous eigenvectors $\alpha, \beta \in \C^4$
of $A$ and $B$.
Then, $(\alpha_i, \beta_i), i = 1, 2, 3, 4$ are linearly independent over $\R$, defining
a lattice of $\Delta$ preserved by $\phi(\lambda)$ and $\phi(\mu)$ 
($\lambda = \log \gamma, \mu = \log \delta$).
And thus we have determined all lattices of $G$.

\begin{remark}
Since $\lambda$ and $\mu$ are linearly independent over $\R$,
we have either $|\gamma| \not=1$ or $|\delta| \not= 1$. And if, for instance,
$|\gamma| \not=1$ and $\gamma \notin \R$,
then $A$ has four distinct eigenvalues
$\gamma^{-1}, \gamma, \overline{\gamma}^{-1}, \overline{\gamma}$.
For the case where both $A$ and $B$ have real eigenvalues
$\gamma^{-1}, \gamma$ and $\delta^{-1}, \delta$ respectively,
take simultaneous non-real eigenvectors $\alpha, \beta \in \C^4$ for them;
then we see that
$(\alpha_i, \beta_i), i = 1, 2, 3, 4$ are linearly independent over $\R$, defining
a lattice $\Delta$ of $\C^2$ preserved by $\phi(\lambda)$ and $\phi(\mu)$.
\end{remark}

\begin{example}
Take $A \in {\rm SL} (4, \Z)$ with four non-real eigenvalues
$\gamma, \gamma^{-1},$ $\overline{\gamma}, \overline{\gamma}^{-1}$;
for instance,
$$A = \left(
\begin{array}[c]{cccc}
0 & 1 & 0 & 0\\
0 & 0 & 1 & 0\\
0 & 0 & 0 & 1\\
-1 & 1 & -3 & 1
\end{array}
\right),
$$
with the characteristic polynomial given by
$${\rm det} (t I - B) = t^4 - t^3 + 3 t^2 - t + 1.$$
For the lattice $\Lambda$ of $\C$ generated by
$\lambda = \log \gamma$ and
$\mu = k \pi \sqrt{-1}\;(k \in~\Z)$,
and the lattice $\Delta$ of $\C^2$ generated by
$(\alpha_i, \beta_i), i = 1, 2, 3, 4$,
we can define a lattice $\Gamma = \Delta \rtimes \Lambda$ of $G$,
where $(\alpha_1, \alpha_2, \alpha_3, \alpha_4),
(\beta_1, \beta_2, \beta_3, \beta_4) \in \C^4$ are eigenvectors
of $A$ with eigenvalue $\gamma, \gamma^{-1}$.
\end{example}

\begin{example}[\cite{N}]
Take $A \in {\rm SL}(2, \Z)$ with two real
eigenvalues $\gamma^{-1}, \gamma$, $\gamma \not= \pm 1$, and
their real eigenvectors 
$(a_1, a_2), (b_1, b_2) \in \R^2$. 
Then, for any $\epsilon \notin \R$ (e.g. $\epsilon = \sqrt{-1}$),
$(a_1, a_2, a_1 \epsilon, a_2 \epsilon)$
and $(b_1, b_2, b_1 \epsilon, b_2 \epsilon)$ are
non-real eigenvectors for 
$A \oplus A \in {\rm SL}(4, \Z)$ with eigenvalues $\gamma^{-1}, \gamma$.
For the lattice $\Lambda$ of $\C$ generated by
$\lambda\;(\lambda = \log \gamma)$ and
$\mu = k \pi \sqrt{-1}\;(k \in~\Z)$,
and the lattice $\Delta$ of $\C^2$ generated by $(a_1, b_1),
(a_2, b_2), (a_1 \epsilon, b_1 \epsilon), (a_2 \epsilon, b_2 \epsilon)$,
we define a lattice $\Gamma = \Delta \rtimes \Lambda$ of $G$.
\end{example}
\end{list}

Let $M = \Gamma \backslash G$ be a three-dimensional  complex
solvmanifold, where $G$ is a simply connected solvable Lie group
with lattice $\Gamma$.
Then, since $G$ is linear algebraic, applying a fundamental theorem of
Winkelmann \cite{Wi}, we have
$${\rm dim}\, H^1(M, {\mathcal O}) = {\rm dim}\, H^1(\g, \C) + {\rm dim}\, W,
$$
where ${\mathcal O}$ denotes the structure sheaf of $M$, $\n$ the nilradical of
$\g$, and $W$ the maximal linear subspace of $[\g, \g]/[\n, \n]$ such that
${\rm Ad}(\xi)$ on $W$ is a real semi-simple linear endomorphism
for any $\xi \in \Gamma$. 
Note that ${\rm dim}\, H^1(\g, \C) = {\rm dim}\,\g - {\rm dim}\,[\g, \g]$, and
${\rm Ad}(\xi)|W$ is diagonalizable over~$\R$.
\smallskip

We can determine $h^1 = {\rm dim}\, H^1(M, {\mathcal O})$ completely
from Winkelmann's formula above and our classification of three-dimensional 
complex solvmanifolds (cf.~\cite{N}):

\begin{list}{}{\topsep=5pt \leftmargin=20pt \itemindent=5pt \parsep=0pt \itemsep=3pt}
\item[ (1)\,] Abelian Type:  ${\rm dim}\, W = 0$, $h^1 = 3$;
\item[ (2)\,] Nilpotent Type: ${\rm dim}\, W = 0$, $h^1 = 2$;
\item[ (3a)] Non-Nilpotent Type with either $\gamma\; {\rm or}\; \delta \notin \R$:
${\rm dim}\, W = 0$, $h^1 = 1$;
\item[ (3b)] Non-Nilpotent Type with $\gamma, \delta \in \R$:
${\rm dim}\, W = 2$, $h^1 = 3$;
\end{list}

We see that complex solvmanifolds in Example 2 are of type (3a),
and those in Example 3 are of type (3b).
\begin{remark}
There seems an error in the construction of a lattice in
the example of  a complex solvmanifold of type (3a)
in the paper \cite{N}.
\end{remark}


\section{Pseudo-K\"ahler structures on complex solvmanifolds}

We recall the definition of pseudo-K\"ahler structure. 
Let $M$ be a symplectic manifold with symplectic form $\omega$. If $M$
admits a complex structure $J$ such that $\omega(JX, JY) = \omega(X, Y)$
for any vector fields $X, Y$ on $M$, we call $(\omega, J)$ a {\em pseudo-K\"ahler}
structure on $M$. For a pseudo-K\"ahler structure $(\omega, J)$, we have a
pseudo-Riemannian structure $g$ defined by
$g(X, Y) = \omega(X, JY)$; if, in addition, $g$
is Riemannian (i.e. positive definite), then we call $(\omega, J)$ a K\"ahler structure
on $M$. Equivalently, a pseudo-K\"ahler (K\"ahler) structure is nothing but
a pseudo-Hermitian (Hermitian) structure with its closed fundamental form $\omega$.

\begin{theorem}
A three-dimensional  complex solvmanifold
admits a pseudo-K\"ahler
structure if and only if it is of type (1), or of type (3b).
\end{theorem}

\begin{proof}
It is known (due to Yamada \cite{Y1}) that a 
 complex solvmanifold of dimension $n$ 
with pseudo-K\"ahler structure must have
$h^1  \ge n$ (actually the equality holds here); in particular,
a  complex solvmanifold of nilpotent type or non-nilpotent type with
either $\gamma\; {\rm or}\; \delta \notin \R$ admits no pseudo-K\"ahler structures.
Therefore, in order to prove the theorem, it is sufficient to show that 
a  complex solvmanifold $\Gamma \backslash G$
of non-nilpotent type with
$\gamma, \delta \in \R$ admits a pseudo-K\"ahler structure.
In Section~3, we observed that we have $\gamma, \delta \in \R$ if and only if
$\Lambda$ is generated by $\lambda = a + k \pi \sqrt{-1}, \mu = b + l \pi \sqrt{-1}$,
where $a, b \in \R\; {\rm and}\; k, l \in \Z$. We can construct a
pseudo-K\"ahler structure $\omega$ on $\Gamma \backslash G$,
as in the paper \cite{Y1}, in the following:
$$ \omega = \sqrt{-1} dx \wedge d \overline{x} + dy \wedge d \overline{z}
+ d \overline{y} \wedge d z,
$$
or using Maurer-Cartan forms (left-invariant 1-forms)
$\omega_1, \omega_2, \omega_3$,on $G$,
$$ \omega = \sqrt{-1} \omega_1 \wedge \overline{\omega_1} +
e^{- 2\, {\rm Im} (x) \sqrt{-1}}\, \omega_2 \wedge \overline{\omega_3} +
e^{2\, {\rm Im} (x) \sqrt{-1}}\, \overline{\omega_2} \wedge \omega_3,
$$
where $\omega_1 = dx, \omega_2 = e^x\, dy, \omega_3 = e^{-x}\, dz$.
\end{proof}
\smallskip

\begin{remark}
We know (\cite{W1}, \cite{H2}) that a complex
solvmanifold admits K\"ahler structures if and only if it is a complex torus.
On the other hand, we know  \cite{DG} that  a complex solvmanifold
admits homogeneous (invariant) pseudo-K\"ahler structures if and only if
it is a complex torus.
Therefore,  a complex solvmanifold, except a complex
torus, admits neither K\"ahler nor homogeneous pseudo-K\"ahler structures.

\end{remark}

\begin{remark}
In the paper \cite{Y2} Yamada showed, applying  Winkelmann's formula,
that a homogeneous complex pseudo-K\"ahler solvmanifold
has the structure of complex torus bundle over a complex torus.
\end{remark}


\section{Left-invariant complex structures on complex solvmanifolds}

Let $G$ denote a complex solvable Lie group of non-nilpotent
type (as defined in section 3), and $\g$ its Lie algebra.
Recall that $\g$ has a basis $X, Y, Z$ over $\C$ with bracket multiplication defined by
$$[X, Y] = -Y,\; [X, Z] = Z,\; [Y, Z] = 0.
\eqno{(1)}$$

Let $\g_\R$ denote the real Lie algebra underlying $\g$. Then,
$\g_\R$ has a basis $X, X', Y, Y', Z, Z'$ over $\R$ with bracket multiplication defined by
$$[X, Y] = -Y,\; [X, Y'] = -Y'\;, [X, Z] = Z,\; [X, Z'] = Z',
\eqno{(2a)}$$
\vspace{-15pt}
$$[X', Y] = -Y',\; [X', Y'] = Y,\; [X', Z] = Z',\; [X', Z'] = -Z,
\eqno{(2b)}$$
and all other brackets are 0.

Let $\g_\C$ denote the complexification of $\g_\R$, that is,
$$\g_\C = \g_\R \oplus \sqrt{-1} \g_\R.
$$

We have the following split short exact sequence:
$$0 \longrightarrow {\mathfrak a}  \stackrel{i}{\longrightarrow} {\g_\C}
\stackrel{r}\longrightarrow
{\mathfrak b} \longrightarrow 0,
$$
where ${\mathfrak a} = [\g_\C, \g_\C]$,
and ${\mathfrak b}$ is the Lie subalgebra of $\g_\C$ generated by $X, X'$ over $\C$.
\smallskip

We now suppose that $\g$ has a left-invariant complex structure $J$ with its
associated complex subalgebra 
$\h$ of $\g_\C$ such that $\g_\C = \h \oplus \overline{\h}$. Then,
${\mathfrak q} = r(\h)$ has the dimension 1 or 2, and
${\mathfrak k} = {\rm ker}(r | \h)$ has the dimension 2 or 1  accordingly.
But we see that
the second case is not possible, and thus we have the following split short exact
sequence:
$$0 \longrightarrow {\mathfrak k}  \stackrel{i}{\longrightarrow} {\h}
\stackrel{r}\longrightarrow
{\mathfrak q} \longrightarrow 0,
$$
where ${\rm dim}\, {\mathfrak k} = 2, {\rm dim}\, {\mathfrak q} = 1$, and
${\mathfrak a} = {\mathfrak k} \oplus \overline{\mathfrak k}$, 
${\mathfrak b} = {\mathfrak q} \oplus \overline{\mathfrak q}$. We can further assume that
${\mathfrak q}$ is generated by $U + \sqrt{-1} U'$ over $\C$, and
${\mathfrak k}$ is generated by $V + \sqrt{-1} V', W + \sqrt{-1} W'$ over $\C$ such that
$$(U, U') = (X, X')\, Q,\; (V, V', W, W') = (Y, Y', Z, Z')\, P,
\eqno{(3)}$$
for some $Q = (q_{i j}) \in {\rm GL}(2, \R)$ and $P = (p_{k l}) \in {\rm GL}(4, \R)$.
And since $\h$ is a subalgebra of $\g_\C$, the following condition must be satisfied:
$$ [U + \sqrt{-1}\, U', V + \sqrt{-1}\, V'] \; = 
2 \alpha\, (V + \sqrt{-1}\, V') +  2 \beta\, (W + \sqrt{-1}\, W'),
\eqno{(4a)}
$$
\vspace{-20pt}
$$[U + \sqrt{-1}\, U', W + \sqrt{-1}\, W']  = 
2 \gamma\, (V + \sqrt{-1}\, V') +  2 \delta\, (W + \sqrt{-1}\, W'),\;\;
\eqno{(4b)}
$$
for some $\alpha, \beta,\gamma, \delta \in \C$.
\medskip

Remark that for the case $Q = I \in {\rm GL}(2, \R)$ and
$P = I \in {\rm GL}(4, \R)$, we have
$$[X + \sqrt{-1}\, X', Y + \sqrt{-1}\, Y']  = 
- 2\, (Y + \sqrt{-1}\, Y'),
$$
\vspace{-20pt}
$$[X + \sqrt{-1}\, X', Z + \sqrt{-1}\, Z']  = \;\;
2\, (Z + \sqrt{-1}\, Z'),
$$
which defines the original complex structure $J_0$ on $G$ (as a complex Lie group)
with its associated complex subalgebra $\h_0$ generated by
$$X + \sqrt{-1}\, X', Y + \sqrt{-1}\, Y', Z + \sqrt{-1}\, Z'$$
 over $\C$.

\begin{lemma}
{\em Let $A = \left(
\begin{array}[c]{cc}
\alpha & \beta\\
\gamma & \delta\\
\end{array}
\right)
\in {\rm GL}(2, \C)$, which satisfies the above equations $(4a), (4b)$.
Then, for ${\mathfrak q}$ and ${\mathfrak k}$ to be Lie subalgebras of 
${\mathfrak g}_{\C}$, $A$ must be conjugate over $\R$ to
$\frac{q_{1 1} + q_{ 2 2}}{2}
\left(
\begin{array}[c]{cc}
-1 & 0\\
0 & 1\\
\end{array}
\right)$
with $q_{1 1} + q_{ 2 2} \not= 0$.}
\end{lemma}

\begin{proof}
For simplicity, we divide the proof into three steps.
\smallskip

\noindent [Step 1]
$P = I$ and $Q = (q_{i j}) \in {\rm GL}(2, \R)$, satisfying the equations $(2a), (2b)$
\smallskip

In this case, we see by calculation that $Q$ is symmetric (i.e. $q_{1 2} = q_{2 1}$),
and we have

$$[U + \sqrt{-1}\, U', Y + \sqrt{-1}\, Y']  = 
-  \frac{q_{1 1} + q_{ 2 2}}{2}\, (Y + \sqrt{-1}\, Y'),
$$
\vspace{-10pt}
$$[U + \sqrt{-1}\, U', Z + \sqrt{-1}\, Z']  = \;\;\;
\frac{q_{1 1} + q_{ 2 2}}{2}\, (Z + \sqrt{-1}\, Z').
$$
\smallskip

\noindent [Step 2]
$Q = I$ and $P = (p_{k l}) \in {\rm GL}(4, \R)$, satisfying the equations $(2a), (2b)$
\smallskip

First, we define a linear automorphism $T \in {\rm Aut}({\mathfrak k})$ by
$$T\, (V, V', W, W') = (V, V', W, W') 
\left(
\begin{array}[c]{cc}
J & O\\
O & J\\
\end{array}
\right),
$$
where $J = \left(
\begin{array}[c]{cc}
0 & 1\\
-1 & 0\\
\end{array}
\right)$,
and a linear endomorphism $S_{X X'} \in {\rm End}({\mathfrak k})$ by
$$S_{X X'} = \frac{1}{2} ({\rm ad}\,X + {\rm ad}\, X' \circ T).
$$
Then,  we have the equation:
$$S_{X X'}\, (Y, Y', Z, Z') = (Y, Y', Z, Z')\, 
\left(
\begin{array}[c]{cc}
-I & O\\
O & I\\
\end{array}
\right),
\eqno {(5)}
$$
and since $\h$ is a subalgebra of $\g_\C$,
the following equation must be also satisfied:
$$S_{X X'}\, (V, V', W, W') = (V, V', W, W')\, 
\left(
\begin{array}[c]{cc}
\alpha I & \beta I\\
\beta I & \delta I\\
\end{array}
\right).
\eqno{(6)}
$$
where  $I \in {\rm GL} (2, \C)$. 
Recall that we have defined $P \in {\rm GL}(4, \R)$ as
$$(V, V', W, W') = (Y, Y', Z, Z')\, P.
\eqno{(7)}
$$
Hence, from the equations $(1), (2), (3)$, we get the equation:
$$\left(
\begin{array}[c]{cc}
-I & O\\
O & I\\
\end{array}
\right)\,
P
= P\,
\left(
\begin{array}[c]{cc}
\alpha I & \beta I\\
\beta I & \delta I\\
\end{array}
\right).
$$
It follows, by simple Linear algebra, that 
$A$
is conjugate over $\R$ to
$\left(
\begin{array}[c]{cc}
-1 & 0\\
0 & 1\\
\end{array}
\right)$.
\smallskip

\noindent [Step 3]
The general case for $P = (p_{k l}) \in {\rm GL}(4, \R), 
Q = (q_{i j}) \in {\rm GL}(2, \R)$, satisfying the equations $(2a), (2b)$
\smallskip

Following the arguments in the Step 1 and 2,
for any $Q \in {\rm GL}(2, \R)$ and $P \in {\rm GL}(4, \R)$,
we see by calculation that $Q$ is symmetric, and 
$A = \left(
\begin{array}[c]{cc}
\alpha & \beta\\
\gamma & \delta\\
\end{array}
\right) \in {\rm GL}(2, \C)$
is conjugate over $\R$ to
$\frac{q_{1 1} + q_{ 2 2}}{2}
\left(
\begin{array}[c]{cc}
-1 & 0\\
0 & 1\\
\end{array}
\right)$.
\smallskip

\noindent This complete the proof of Lemma 2.
\end{proof}
\medskip

We see, from Lemma 2, that there exists a complex automorphism of
Lie algebras $\Phi:\g_\C \rightarrow \g_\C$ such that
$\Phi \circ \tau_0 = \tau \circ \Phi$ and $\Phi (\h_0) = \h$,
where $\tau$ and $\tau_0$ are the conjugations with respect to $J$ and $J_0$
respectively. In fact, we have an equivalence of short exact sequences

$$\begin{array}{ccccccccc}
0 & \longrightarrow & {\mathfrak k}_0  & \stackrel{i}{\longrightarrow} & {\h}_0
& \stackrel{r}\longrightarrow
& {\mathfrak q}_0 & \longrightarrow & 0\\
&  & \vcenter{\rlap{$\scriptstyle k$}}\: \downarrow & &
\vcenter{\rlap{$\scriptstyle h$}}\: \downarrow 
& &  \vcenter{\rlap{$\scriptstyle q$}}\: \downarrow 
& &\\
0 & \longrightarrow & {\mathfrak k}  & \stackrel{i}{\longrightarrow} & {\h}
& \stackrel{r}\longrightarrow
& {\mathfrak q} & \longrightarrow & 0
\end{array}
$$
satisfying ${\rm ad}(q(u)) \circ k = k \circ {\rm ad} (u)\;
(u = U + \sqrt{-1}\, U')$, 
which extends to  a complex automorphism of
Lie algebras $\Phi:\g_\C \rightarrow \g_\C$ such that
$\Phi \circ \tau_0 = \tau \circ \Phi$ and $\Phi (\h_0) = \h$.
To be more precise,  for  $K  \in {\rm GL}(2, \R)$ such that
$$K^{-1} A K = \frac{q_{1 1} + q_{ 2 2}}{2}
\left(
\begin{array}[c]{cc}
-1 & 0\\
0 & 1\\
\end{array}
\right),$$
$k$ is a linear map defined by $K$, and $q$ is a scalar multiplication by 
$\frac{2}{q_{1 1} + q_{ 2 2}}$.

\smallskip

For the original complex solvable Lie group $(G, J_0)$ with
its associated complex subalgebra $\h_0$, the complex subgroup $H_0$
of $G_\C$ corresponding to $\h_0$  is closed, simply connected and
$H_0 \backslash G_\C$ is biholomorphic to $\C^3$. 
We have $\Gamma = G \cap H_0 = \{1\}$, 
and the canonical map $ g_0 = q_0 \circ i$
$$(G, J_0) \stackrel{i}{\hookrightarrow}  G_\C \stackrel{q_0}{\rightarrow}
H_0 \backslash G_\C
$$
is a biholomorphic map. The complex automorphism of Lie algebras
$\Phi$] induces a complex automorphism of Lie groups
$\Psi: G_\C \rightarrow G_\C$ such that $q \circ \Psi= \tilde{\Psi} \circ q_0$,
which send $H_0$ to $H$ biholomorphically. 
$$\begin{array}{ccccc}
(G, J_0) & \stackrel{i}{\hookrightarrow} & G_\C & \stackrel{q_0}{\rightarrow} &
H_0 \backslash G_\C \\
 &  & \vcenter{\rlap{$\scriptstyle \Psi$}}\: \downarrow & &
\vcenter{\rlap{$\scriptstyle \tilde{\Psi}$}}\: \downarrow \\
(G, J) & \stackrel{i}{\hookrightarrow} & G_\C & \stackrel{q}{\rightarrow} &
H \backslash G_\C 
\end{array}
$$
Hence, the canonical map $g = q \circ i$
is also a biholomorphic map.
We have thus shown
\smallskip

\begin{proposition}
Let $G$ be  a three-dimensional simply connected
complex solvable Lie group of non-nilpotent type. Then,
any left-invariant complex structure on $G$ is biholomorphic to $\C^3$.
\end{proposition}

We can also show  that any left-invariant complex structures on
three-dimensional complex solvable Lie groups of abelian type or nilpotent type 
are biholomorphic to $\C^3$. The proof is almost the same as for the case of
non-nilpotent type in Theorem 3. We have thus shown 

\begin{theorem}
Any left-invariant complex structure on a three-dimensional simply connected
complex solvable Lie group is biholomorphic to $\C^3$.
\end{theorem}

We know (due to Kodaira \cite{N}) that 
among small deformations of a three-dimensional  complex
solvmanifold of type (3b) there exists a continuous
family of complex structures whose universal coverings are not Stein.
We have thus obtained 

\begin{theorem}
There exists a continuous family of non-left-invariant complex structures on
a three-dimensional  complex solvmanifold of type (3b)
\end{theorem}

\begin{remark}
Let $M = \Gamma \backslash G$ be an Iwasawa manifold, a three-dimensional  complex
nilmanifold (complex solvmanifold of type (2)), where $\Gamma$  is a lattice of
a simply connected nilpotent Lie group $G$.
It is known (\cite{N}, \cite{S}) that the moduli space of all left-invariant
complex structures on an Iwasawa manifold has the dimension 6,
while small deformations (Kuranishi space) of the Iwasawa manifold
also has the dimension 6; and all of their universal coverings are biholomorphic
to $\C^3$. It follows that small deformations of the Iwasawa manifold are all
left-invariant. 
\end{remark}


\begin{flushleft}
Department of Mathematics\\
Faculty of Education\\
Niigata University, Niigata\\
JAPAN\\
\vskip3pt
e-mail: hasegawa@ed.niigata-u.ac.jp\\
\end{flushleft}

\end{document}